\RequirePackage{fix-cm}
\documentclass[smallcondensed]{svjour3}

\smartqed
\usepackage{graphicx}
\usepackage{latexsym}
\usepackage[utf8]{inputenc}
\usepackage{amssymb}
\usepackage{mathtools}
\usepackage{mathrsfs}
\usepackage{emptypage}

\newtheorem{pr}[theorem]{Proposition}
\newtheorem{df}[theorem]{Definition}

\newtheorem{co}[theorem]{Corollary}

\newtheorem{lm}[theorem]{Lemma}

\begin{document}

\title{On the existence of Ulrich bundles on blown-up varieties at a point}

\author{Saverio Andrea Secci}

\institute{S. A. Secci}

\date{}


{\center{\textbf{\Large{On the existence of Ulrich bundles on blown-up varieties at a point}}}}
\bigskip

{\center{\textbf{\large{Saverio Andrea Secci}}}}
\vspace{3cm}

\begin{abstract}
The objective is to show the construction of an Ulrich vector bundle on the blowing-up $\widetilde X$ of a nonsingular projective variety $X$ at a closed point, where the original variety is embedded by a very ample divisor $H$ and carries an Ulrich vector bundle. In order to achieve this result, we aim to find a suitable very ample divisor on $\widetilde X$, which is dependent on $H$.
At the end, we take into consideration some applications to surfaces with regards to minimal models and their Kodaira dimension.
\end{abstract}

\section{Introduction}

Ulrich vector bundles made their first appearance during the 80's in commutative algebra, thanks to Ulrich \cite{U}, who studied them under the name “maximally generated maximal Cohen–Macaulay modules”.
In the latest years, due to their connection with other notions, such as determinantal hypersurfaces and Chow forms (see \cite{ES} and \cite{B1}), Ulrich bundles have obtained considerable attention.
Although there exist various characterisations, we define them as follows:

\begin{df}
Let $X$ be a smooth projective variety of dimension $n$, and let $H$ be a very ample divisor on $X$. Let $\mathcal E$ be a vector bundle on $X$ such that, for all $1 \leq p \leq n$, $\mathcal E(-pH)$ has vanishing cohomology. Such $\mathcal E$ is said to be an Ulrich vector bundle for $(X, H)$.
\end{df}

The question of existence of Ulrich vector bundles, despite many results being already proven, is still open.
In the case of curves, finding an Ulrich line bundle $\mathcal E$ on a smooth projective curve $C$ of genus $g$ is equivalent to finding a divisor $D$ of degree $g-1$ with no global sections, and indeed they are in a one-to-one correspondence. As a matter of fact, there are divisors satisfying such properties: if we consider the map $\varphi : C^{g-1} \to {\rm Pic}^{g-1}(C)$, which sends $g-1$ points $(p_1, \dots, p_{g-1})$ to the line bundle $\mathcal O_C(p_1+ \dots +p_{g-1})$, then it is known that $\varphi$ is not surjective, and therefore we find the requested divisors in ${\rm Pic}^{g-1}(C) \backslash {\rm Im}( \varphi)$. In fact, if we take a line bundle $\mathcal L := \mathcal O_C(D)\in {\rm Pic}^{g-1}(C) \backslash {\rm Im}( \varphi)$, we have that $h^0(\mathcal L)=0$ and that $\deg (\mathcal L)=g-1$, therefore $h^1(\mathcal L)=0$ by Riemann-Roch theorem. Hence the line bundle $\mathcal E=\mathcal L(1)$ is Ulrich for $(C, \mathcal O_C(1))$.

More intricate is, instead, the case of dimension two. We already know that several classes of surfaces admit an Ulrich vector bundle, such as abelian surfaces \cite{B2}, Enriques surfaces \cite{B3,C1}, K3 surfaces \cite{F}, bielliptic surfaces \cite{B4}, many geometrically ruled surfaces \cite{ACM}, surfaces with $p_g=0$, $q=1$ \cite{C2}, surfaces with $q=0$ and arbitrary $p_g$, which satisfy few technical conditions \cite{C3}, hypersurfaces and complete intersections \cite{HUB}, surfaces with $q=1$ and surfaces with $q \geq 2$ and minimal model $S_0$ such that ${\rm Pic}(S_0) \cong \mathbb Z$ \cite{Lo}.

Our intent is to analyse the problem of existence with regards to blown-up varieties: how does the blowing-up of a nonsingular projective variety $X$ at a closed point behave with respect to the existence of an Ulrich vector bundle for $(X, H)$? Taking into consideration \cite[Thm.0.1]{K}, we provide the following explicit version:

\begin{theorem}
\label{MyThm}
Let $H$ be a very ample divisor on a nonsingular projective variety $X$. Let $P \in X$ be a closed point, and let $\pi:\widetilde X \to X$ be the blowing-up of $X$ at $P$. If there exists an Ulrich vector bundle for $(X, 2H)$, then there exists an Ulrich vector bundle for $(\widetilde X, 2\pi^*H -E)$, where $E=\pi ^{-1}(\{P\})$ is the exceptional divisor.
\end{theorem}

Note that in \cite[Thm.0.1]{K} the hypothesis on $H$ is that $\pi^*H-E$ is very ample, which does not hold in general. It holds, for example, when $H$ is sufficiently ample. Theorem \ref{MyThm}, on the other hand, provides an Ulrich vector bundle for $(\widetilde X, 2\pi^*H -E)$ regardless of the properties of $H$ since, if $H$ is very ample, then $2\pi^*H -E$ is always very ample, as proved in Corollary \ref{co6}. Let us note that, while a direct proof for Theorem \ref{MyThm} is provided, one can observe that Theorem \ref{MyThm} is also a consequence of both \cite[Thm.0.1]{K} and Corollary \ref{co6}.

Moreover, although \cite{K} contains a gap - which is underlined and corrected in \cite{CK} - it does not affect \cite[Thm.0.1]{K}.\\

\textbf{Remark}
Theorem \ref{MyThm} is also true if we suppose the existence of an Ulrich vector bundle for $(X, H)$.
This follows by \cite[Proposition 5.4]{ES} (see also \cite[Section 3]{B4}), which states that if there exists an Ulrich vector bundle for $(X, H)$, then there exists an Ulrich vector bundle for $(X, 2H)$.\\

Furthermore, Theorem \ref{MyThm} yields interesting applications to surfaces. Since every smooth surface is obtained by a finite number of blow-ups of its minimal surface at closed points (see \cite[Chapter II]{B}), we deduce that in order to investigate the question of existence of Ulrich vector bundles in the case of dimension two we can focus on studying minimal models.

\begin{co}
\label{MyCor}
If every minimal model carries an Ulrich vector bundle, then for every nonsingular projective surface $S$ there exists a very ample divisor for which $S$ admits an Ulrich vector bundle.
\end{co}

One way for approaching the problem is to consider surfaces based on their Kodaira dimension, a birational invariant, and we are going to analyse surfaces with Kodaira dimension $\leq 0$.
Thanks to other results and some computation showed in Section 3, we are going to see that for every minimal model $S_0$ with Kodaira dimension $\leq 0$ there exists a very ample divisor for which $S_0$ carries an Ulrich vector bundle.

Clearly, the very ample divisors arising from this method become more and more high as one blows-up the underlying variety.
That is, let $\pi_i:X_i \to X$ the blow-up of $X$ at the points $P_1, \dots P_i$. The very ample divisor on $X_1$ for which there exists an Ulrich bundle is $2\pi_1^*H -E_1$; on $X_2$ we have $4\pi_2^*H- 2E_1-E_2$; on $X_i$ we have $2^i\pi_i^*H - \sum_{j=1}^i 2^{i-j}E_j$.
In general there may be many very ample divisors for which existence of Ulrich vector bundles is still unknown.

If not specified otherwise, a $scheme$ is a separated scheme of finite type over an algebraically closed field $K$. A $variety$ is an integral scheme.

\section{Ulrich vector bundle on blowing-ups}
Our first goal is to determine a suitable very ample divisor on $\widetilde X$.

The following result is taken from \cite[Thm.2.1]{BS}. Note that, although \cite{BS} works on the field of complex numbers, one can observe that \cite[Thm.2.1]{BS} does not strictly require that the ground field is $\mathbb C$ to be true, therefore we can apply it even in our more general case.

\begin{theorem}
\label{2.1BelSom}
Let $\mathcal L$ be a very ample line bundle on a nonsingular projective variety $X$, and let $Y \subseteq X$ be a closed subscheme corresponding to a sheaf of ideals $\mathscr I$ on $X$.
Let $\pi:\widetilde X \to X$ be the blowing-up of $X$ with respect to $\mathscr I$, and let $E=\pi ^{-1}(Y)$ be the exceptional divisor.
Assume that $\mathscr I \otimes \mathcal L^{\otimes t}$ is generated by global sections for some positive integer t. Then $t'\pi ^* \mathcal L -E$ is very ample for $t' \geq t+1$.
\end{theorem}

\begin{pr}
\label{IxL.GG}
Let X be a scheme, and let $\mathcal L$ be a very ample line bundle on X corresponding to a closed immersion $i:X \to \mathbb P^N$. Let $P \in X$ be a closed point. Then $\mathscr I_{\{P\} /X} \otimes \mathcal L$ is globally generated.
\end{pr}
\textit{Proof:}
At first we observe that $\mathscr I_{\{P\} / \mathbb P^N}(1)$ is globally generated: this follows by the fact that the ideal $I_{\{P\}}=\Gamma_*(\mathscr I_{\{P\} / \mathbb P^N})$ of a point $P \in \mathbb P^N$ is generated by $N$ independent hyperplanes.

Let us now consider the surjective morphism
\[\mathscr I_{\{P\} / \mathbb P^N}(1) \to i_* \mathscr I_{\{P\} /X} \otimes \mathcal O_{\mathbb P^N}(1).\]
Since $\mathscr I_{\{P\} / \mathbb P^N}(1)$ is globally generated, we have that $i_* \mathscr I_{\{P\} /X} \otimes \mathcal O_{\mathbb P^N}(1)$ is globally generated.
Furthermore we observe that by the projection formula \cite[Ex.II.5.1(d)]{H} we have $i_* \mathscr I_{\{P\} /X} \otimes \mathcal O_{\mathbb P^N}(1) \cong i_*(\mathscr I_{\{P\} /X} \otimes i^* \mathcal O_{\mathbb P^N}(1))$.\\
This implies that $\mathscr I_{\{P\} /X} \otimes i^* \mathcal O_{\mathbb P^N}(1)=\mathscr I_{\{P\} /X} \otimes \mathcal L$ is globally generated.
\qed

\begin{co}
\label{co6}
Let $H$ be a very ample divisor on a nonsingular projective variety $X$, and let $P \in X$ be a closed point.
Let $\pi:\widetilde X \to X$ be the blowing-up of $X$ at $P$, and let $E=\pi ^{-1}(\{P\})$ be the exceptional divisor. Then $t\pi^*H -E$ is a very ample divisor on $\widetilde X$ for all $t \geq 2$.
\end{co}
\textit{Proof:}
Let $\mathscr I$ be the sheaf of ideals corresponding to $\{P\}$. By Proposition \ref{IxL.GG}, $\mathscr I \otimes \mathcal O_X(H)$ is globally generated. Then, by applying Theorem \ref{2.1BelSom}, we achieve our goal.
\qed 

The following result is taken from \cite[Lemma 1.4]{BEL}.

\begin{lm}
\label{1.4BEL}
Let $X$ be a nonsingular projective variety, and let $\pi:\widetilde X \to X$ be the blowing-up of $X$ at a closed point $P \in X$. Then, if $E$ is the exceptional divisor, for any locally free sheaf $\mathcal F$ on $X$ we have
\[H^i(\widetilde X, \pi^* \mathcal F \otimes \mathcal O_{\widetilde X}(tE))=H^i(X,\mathcal F)\]
for all $i \geq 0$ and for $0 \leq t \leq \dim X -1$.
\end{lm}

\textit{Proof of Theorem \ref{MyThm}:}
Let us fix $n=\dim X$, and let $\mathcal F$ be an Ulrich vector bundle for $(X, 2H)$.
We aim to prove that $\mathcal G :=\pi^* \mathcal F \otimes \mathcal O_{\widetilde X}(-E)$ is an Ulrich vector bundle for $(\widetilde X, 2\pi^*H -E)$.
Indeed:
\[\mathcal G(-p(2\pi^*H -E))=\pi^*\mathcal F(-2pH) \otimes \mathcal O_{\widetilde X}((p-1)E),\]
and by Lemma \ref{1.4BEL} we have that
\[H^i(\widetilde X, \pi^*\mathcal F(-2pH) \otimes \mathcal O_{\widetilde X}((p-1)E))=H^i(X, \mathcal F(-2pH))\]
for all $i \geq 0$ and for $0 \leq p-1 \leq n-1$.
Therefore, since $H^i(X, \mathcal F(-2pH))=0$ for all $i \geq 0$ and for $1 \leq p \leq n$, and since $\dim X=\dim{\widetilde X}$, we achieve the expected result.
\qed

\section{Application to surfaces}
Let us now assume that the ground field $K$ is an algebraically closed field of characteristic zero. We are going to look at some application of Corollary \ref{MyCor}.

We first observe that known results give the existence of an Ulrich vector bundle, for a suitable very ample divisor $H$, on any minimal surface of Kodaira dimension $\le 0$.
In fact, by \cite[Thm.VIII.2]{B}, we know that every minimal surface with Kodaira dimension 0 can be one of the following: abelian surfaces, Enriques surfaces, K3 surfaces, and bielliptic surfaces. As it was pointed out in Section 1, we know that an Ulrich vector bundle for a suitable very ample divisor has been found for each of the four cases above.

As for surfaces with Kodaira dimension $-\infty$, by \cite[Prop.III.21, Thm.VI.17]{B} we deduce that they correspond to ruled surfaces.
Then, by \cite[Thm.III.10, Thm.V.10]{B} we have that the minimal models for ruled surfaces are either $\mathbb P^2$ (which carries the Ulrich bundle $\mathcal O_{\mathbb P^2}$ with respect to the very ample line bundle $\mathcal O_{\mathbb P^2}(1)$) or the geometrically ruled surfaces over a nonsingular projective curve $C$, that is the projective bundles $\mathbf P_C(\mathcal E)$, where $\mathcal E$ is a rank 2 vector bundle over $C$.

Without loss of generality we can assume that the bundle $\mathcal E$ is normalised, and let $e=-\deg(\mathcal E)$ be the invariant of $ \mathbf P_C(\mathcal E)$.
Let $g$ be the genus of $C$ and $B$ be a divisor over $C$ such that
\[\deg B > \max \{2g, e + 2g\},\]
and let $F$ be the divisor over $C$ of degree $-e$ associated to the line bundle $\wedge^2 \mathcal E$. Then
\[\deg(B)>2g, \deg(B+F)>2g,\]
therefore $B$ and $B+F$ are very ample by \cite[Cor.IV.3.2(b)]{H}.
Furthermore for all $P \in C$ we have that
\[\deg(B-P)>2g-1, \deg(B+F-P)>2g-1,\]
so
\[H^1(\mathcal O_C(B-P))=0, H^1(\mathcal O_C(B+F-P))=0.\]
Therefore we have that $H=C_0+Bf$ is very ample by \cite[Ex.V.2.11(b)]{H}, where $C_0 \subset \mathbf P_C(\mathcal E)$ is the section associated to the line bundle $\mathcal O_{P_C(\mathcal E)}(1)$ and $f$ is a fibre.
Since $H.f=1$, by the closed embedding $\mathbf P_C(E) \to \mathbb P^{m-1}$, $m=h^0(\mathcal O_{P_C(\mathcal E)}(H))$, given by the linear system $|H|$, the fibres $f$ are sent into lines.
Eventually, by \cite[Prop.4.1(ii)]{B4}, we conclude that $\mathbf P_C(\mathcal E)$ carries an Ulrich bundle.\\

Those previous observations give us the following result:

\begin{co}
Let $S$ be a nonsingular projective surface with Kodaira dimension $k(S) \leq 0$. Then there exists a very ample divisor $H$ such that $S$ admits an Ulrich vector bundle for $(S,H)$.
\end{co}

\section*{Conflict of interest statement}
On behalf of all authors, the corresponding author states that there is no conflict of interest.

\end{document}